\documentstyle[11pt]{article}
\renewcommand{\baselinestretch}{1.3}

\newtheorem {thm}{Theorem}[section]
\newtheorem {prop}[thm]{Proposition}
\newtheorem {lem}[thm]{Lemma}
\newtheorem {cor}[thm]{Corollary}

\newcommand{\bea}{\begin{eqnarray}}
\newcommand{\ba}{\begin{array}}
\newcommand{\bean}{\begin{eqnarray*}}
\newcommand{\ea}{\end{array}}
\newcommand{\eea}{\end{eqnarray}}
\newcommand{\eean}{\end{eqnarray*}}
\newcommand{\be}{\begin{equation}}
\newcommand{\ee}{\end{equation}}

\newcommand{\lab}{\label}

\def\Cox{\hfill \Box}

\def\sif{\sigma \mbox{\rm -field}}

\def\E{{\bf{E}}}
\def\P{{\bf{P}}}
\def\N{{\bf N}}
\def\C{{\cal{C}}}
\def\Z{{\bf{Z}}}
\def\F{{\cal{F}}}

\def\Ga{\Gamma}
\def\B{{\cal{B}}}
\def\|{\, | \, }
\def\v0{{\bf 0}}
\def\one{{\bf 1}}
\def\0{\hat{0}}
\def\1{\hat{1}}
\def\al{\alpha}
\def\la{\lambda}
\def\th{\theta}
\def\pa{\Upsilon}
\def\pas{\Upsilon}
\def\wh{\widehat}
\def\T{{\bf T}}
\def\phi{\varphi}
\def\pr{{\sf PRE}}
\begin{document}
\begin{center}
{\Large \bf Unpredictable Paths and Percolation} 
\end{center}
\vspace{2ex}
\begin{center}
{\sc Itai Benjamini}\,\footnote{Mathematics Department, Weizmann Institute,
Rehovot 76100, Israel.}, \, 
{\sc Robin Pemantle}\,\footnote{Research supported in part by 
National Science Foundation grant \# DMS-9300191, by a Sloan Foundation
Fellowship, and by a Presidential Faculty Fellowship.}$^,$\footnote{Department 
of Mathematics, University of Wisconsin, Madison, WI 53706 \, .},
\,  and \, 
{\sc Yuval Peres}\,\footnote {Research partially supported by NSF grant
\# DMS-9404391}$^,$\footnote{Mathematics Institute, The Hebrew University, 
Jerusalem, Israel
and Department of Statistics,  University of California, Berkeley, CA.}
\end{center}
\vspace{5ex}
\begin{center}
{\bf Abstract}
\end{center}
We construct a nearest-neighbor process $\{S_n\}$ on $\Z$ that is less predictable than 
simple random walk,  in the sense that given the process until time $n$,
the conditional probability that $S_{n+k}=x$ is uniformly bounded by $Ck^{-\alpha}$
for some $\alpha>1/2$. From this process, we obtain a probability 
measure $\mu$ on oriented paths in $\Z^3$ such that the number of intersections of two paths
chosen independently according to $\mu$, has an exponential tail.
(For $d \geq 4$, the uniform measure on oriented paths from the origin
 in $\Z^d$ has this property.)
We show that on  {\em any} graph where such a measure on paths exists,
 oriented percolation clusters are transient
if the retention parameter $p$ is close enough to 1. This yields an extension of a theorem of Grimmett, Kesten and Zhang, who proved that supercritical
 percolation clusters in $\Z^d$ are transient for all $d \geq 3$.
\vfill



 
\newpage

\section{Introduction}
An oriented path from the origin in the lattice $\Z^d$ is determined by a sequence
of vertices $\{y_n\}_{n \ge 0}$ where $y_0=0$ and for each 
$n \ge 1$, the increment $y_n-y_{n-1}$ is one of the $d$ standard basis
vectors. When these increments are chosen independently and uniformly
among the $d$ possibilities, we refer to the resulting random path as a 
{\em uniform random oriented path}.
For  $d \geq 4$, the number of intersections
of two (independently chosen) uniform random oriented paths in $\Z^d$, 
has an exponentially decaying tail. Cox and Durrett (1983) used this fact  
to obtain upper bounds on the critical probability 
$p_c^{\mbox{\rm \small  or }}$ for oriented percolation.
 (They attribute this idea to  H. Kesten.)

In $\Z^3$, however, two independently chosen uniform random oriented paths
 have infinitely many intersections a.s. Perhaps surprisingly, there is a different
measure on oriented paths in $\Z^3$, with exponential tail for the intersection number 
(see Theorem \ref{thm:3d} below). The usefulness of such a measure goes beyond
estimates for $p_c^{\mbox{\rm \small  or}}$, since on {\em any} graph, its existence
implies  that for  $p$ close enough to $1$, a.s.\ some 
infinite  cluster  for oriented percolation is  transient
for simple random walk (see Proposition \ref{prop:trans} below). 
In particular, for sufficiently large $p$, oriented clusters are transient in 
$\Z^d$ for all $d \geq 3$.  
This extends a  theorem of  Grimmett, Kesten and Zhang (1993),
who established transience of
 simple random walk on the infinite cluster of ordinary percolation
in $\Z^d$, $d \geq 3$
(They obtain transience for {\em all} $p>p_c$ but in $\Z^d$
this can be reduced to the case of large $p$ by renormalization,
see Section 2).

We construct the required measure in three dimensions 
   from certain   nearest-neighbor stochastic processes on $\Z$ 
which are ``less predictable than simple random walk''.

\noindent{\bf Definition.}\,
For a sequence of random variables $S=\{S_n\}_{n \geq 0}$ taking values in a 
countable set $V$, we define its {\bf predictability profile} $\{\pr_S (k)\}_{k \geq 1}$ 
by
\be \lab{eq:defpre}
\pr_S (k)= \sup \P[S_{n+k} =x \mid S_0, \ldots , S_n]  \, ,
\ee
where the supremum is over all $x \in V$, all $n \geq 0$ and all histories
$S_0, \ldots , S_n$. 

Thus $\pr_S(k)$ is the maximal chance of guessing
$S$ correctly $k$ steps into the future, given the past of $S$.
Clearly, the predictability profile of simple random walk on $\Z \,$
is asymptotic to  $ C k^{-1/2}$ for some $C>0$.
\begin{thm} \label{thm:smooth}
\begin{description}
\item{{\bf (a)}} For any $\alpha<1$ there exists an integer-valued
stochastic process $\{S_n\}_{n \geq 0}$
such that\/  $|S_n-S_{n-1}|=1$ a.s.\  for all $\, n \geq 1$ and
\be \lab{eq:smooth}
\pr_S(k) \leq C_{\alpha} k^{-\al} \quad \mbox{ \rm for some } C_{\al} < \infty,
  \mbox{ \rm for all } k \geq 1 \, .
\ee
\item{{\bf (b)}} Furthermore, there exists such a process where
 the  $\pm 1$ valued  increments $\, \{S_n-S_{n-1}\}$ 
form a  {\bf stationary ergodic} process.
\end{description}
\end{thm}
Part (b) is not needed for the applications in this paper,
and is included because such processes may have independent interest.
The approach that naturally suggests itself
to obtain processes with a low predictability profile, is to 
use a discretization of fractional Brownian motion;
however, we could not turn this idea into a rigorous construction.
Instead, we construct the processes described in  Theorem \ref{thm:smooth} 
    from a variant of the Ising model
on a regular tree, by summing the spins along the boundary of the tree 
(see Section 4). This may be a case of the principle: ``when you have a hammer,
everything looks like a nail'', and we would be interested to
see alternative constructions.

\noindent{\bf Definitions.}
\vspace{-.1in}
\begin{enumerate}
\item
Let $G=(V_G,E_G)$ be an infinite directed graph with all vertices of
finite degree and let $v_0 \in V_G$. Denote by 
$\pa= \pa (G,v_0)$ the collection of infinite directed paths in $G$ which 
emanate from $v_0$ and tend to infinity (i.e., the paths in $\pa$ visit any vertex 
at most finitely many times). 
 The set $\pa(G,v_0)$, viewed as a subset of
$E_G^\N$, is a Borel set in the product topology. 
\item 
Let $0 < \th <1$. A Borel probability measure $\mu$ on $\pa(G,v_0)$ 
has   {\bf Exponential intersection tails} with parameter $\th$ 
(in short, EIT($\th$)) if there exists $C$ such that
\be \lab{eq:eit} 
\mu \times \mu \Big\{(\varphi,\psi): | \varphi \cap \psi | \geq n \Big\}
\leq C \th^n
\ee
for all $n$, where  $| \varphi \cap \psi |$ is the number of edges in the 
intersection of $\varphi$ and $\psi$.
\item 
If such a measure $\mu$ exists for some  basepoint $v_0$ and some $\th<1$, then
we say that $G$ {\em admits random paths with\/} EIT($\th$).
Analogous definitions apply to undirected graphs.
\item {\bf Oriented percolation}  with parameter $p \in (0,1)$ on the directed graph $G$
is the process where each edge of $G$  is independently declared
 {\em open} with probability $p$ and {\em closed} with probability $1-p$.
The union of all directed open paths emanating from $v$ will be called the 
{\bf oriented open cluster} of $v$ and denoted $\C(v)$. 
\item
A subgraph $\Lambda$ of $G$ is called {\bf transient}
if when the orientations on the edges are ignored,
$\Lambda$ is connected and simple random walk on it is a transient Markov
chain. As explained in Doyle and Snell (1984), the latter property 
is equivalent to
finiteness of the effective resistance from a vertex of $\Lambda$
to infinity, when each edge of $\Lambda$ is endowed with a unit resistor.
\end{enumerate}

\begin{prop} \label{prop:trans}
Suppose a directed graph $G$ admits random paths with EIT($\th$).
Consider oriented percolation on $G$ with parameter $p$.
If $p>\th$ then with probability 1 there is a vertex $v$ in $G$
such that the directed open cluster $\C(v)$ is transient.
\end{prop}
{\sc Remark:} The proof of the proposition, given in Section 2,
also applies to site percolation. If the graph $G$ is {\em undirected} and admits
random (undirected) paths with exponential intersection tails, 
then the same proof shows that  for $p$ close enough to 1, a.s.\ 
some infinite cluster of ordinary percolation on $G$ is transient.  

Recall that a path $\{\Gamma_n\}$ in $\Z^d$ is called {\em oriented} if each 
 increment $\Ga_{n+1}-\Ga_n$ is one of the $d$ standard basis vectors.
 The difference of two independent, uniformly chosen, oriented paths
in $\Z^d$ is a random walk with increments generating the $d-1$
dimensional hyperplane
$\{\sum_{i=1}^d x_i =0 \}$. For $d \geq 4$, this random walk is transient;
let $\th_d<1$ denote  its return probability to the origin. As noted by 
Cox and Durrett (1983), it follows that
 the uniform measure on oriented paths in $\Z^d$ has EIT($\th_d$).
Clearly, a different approach is needed for $d=3$.

\begin{thm} \label{thm:3d}
 There exists a measure on oriented paths from the origin in
 $\Z^3$ that has  exponential intersection tails.
\end{thm}

The rest of this paper is organized as follows. In Section 2 we prove
Proposition \ref{prop:trans} by constructing
a flow of finite energy on the percolation cluster. For $d \geq 4$, this yields a ``soft''
proof that if the parameter $p$ is close enough to 1, 
then oriented percolation clusters in $\Z^d$ are transient with positive probability;
 for $d = 3$, Theorem \ref{thm:3d}
is needed to obtain the same conclusion. We also explain there how transience
of  ordinary percolation clusters in $\Z^d$ for $d \geq 3$ and $\em all\/$ $p>p_c$
can be reduced to transience of these clusters for $p$ close to 1.
In Section 3 we relate the predictability profile of a random process to the tail of 
its collision number with a fixed sequence, and establish Theorem \ref{thm:3d}. 
Theorem \ref{thm:smooth} is proved in Section 4. The main 
ingredient in the proof is an
estimate on the distribution of the population vector in a certain 
two-type branching
process, when this vector is projected in a non-principal eigendirection.
Section 5 contains auxiliary remarks and  problems.
After this paper was submitted, the methods introduced here
were refined and extended by several different authors,
and some of the questions we raised were solved.
The paper concludes with a brief survey of these recent developments.

\section{Exponential intersection tails imply transience of clusters}
To show that an infinite connected graph $\Lambda$
is transient, it suffices to construct a nonzero
 {\bf flow} $f$ on $\Lambda$, with a single source at $v_0$,
such that the {\bf energy} $\sum_{e \in E_\Lambda} f(e)^2$ of $f$ is finite
(see  Doyle and Snell (1984)).

{\sc Proof of Proposition \ref{prop:trans}}: \,
The hypothesis means that there is some vertex $v_0$ and a 
probability measure $\mu$ on $\, \pa= \pa(G, v_0)$ that has EIT($\th$).
 We first assume that 
\be \lab{eq:assump}
\mbox { \it The paths in the closed support of $\: \mu$
 are self-avoiding and tend to infinity uniformly.}
\ee A path is  {\em self-avoiding\/} if it  never revisits a vertex; 
the second part of the assumption means that there is a sequence $r_N \to \infty$
such that for all $N \geq 1$ and all paths $\phi$ in the support of $\mu$,
 the endpoint of $\phi_N$ is not in $B(v_0,r_N)$, where 
 $B(v_0, r)$ denotes the ball of radius $r$ centered at $v_0$ in the usual
graphical distance.
The assumption  (\ref{eq:assump}) certainly  holds in our main application, 
where $\mu$ is supported on
oriented paths in $G=\Z^d$; at the end of the proof
we show how to remove the assumption (\ref{eq:assump}).

For $N \geq 1$ and any infinite path $\varphi \in \pa(G, v_0)$, denote by $\varphi_N$
the finite path consisting of the first $N$ edges of $\varphi$.
Consider the random variable
\be \lab{eq:Z}
Z_N= \int_{\pas} p^{-N} 
\one_{\{\varphi_N \mbox{ \rm is open}\}}  \, d\mu(\varphi) \,.
\ee
Except for the normalization factor $p^{-N}$,
this is the $\mu$-measure of the 
paths that stay in the open cluster of $v_0$ for $N$ steps.

Since each edge is open  with probability $p$ 
(independently of other edges),
$\E(Z_N)=1$.
The second moment of $Z_N$ is
\bea \lab{eq:2nd}
\E(Z_N^2)&=& \E \int_{ \pas} \!\! \int_{ \pas} p^{-2N} 
\one_{\{\varphi_N \mbox{ \rm \small and } 
\psi_N \mbox{ \rm are open}\}}  \, d\mu(\varphi) \, d\mu(\psi) \\[1ex] \nonumber
 &\le &\int_{ \pas} \!\! \int_{ \pas} p^{-|\varphi \cap \psi|}
 \, d\mu(\varphi) \, d\mu(\psi) \,.
\eea
By (\ref{eq:eit}), the last integral is at most 
$
\sum_{n=0}^\infty C \th^n p^{-n}=  C/(1-\th p^{-1}) .
$

 By Cauchy-Schwarz,
\be
\P[|\C(v_0)| \geq N] \geq \P[Z_N>0] \geq {\E(Z_N)^2 \over \E(Z_N^2)}
\geq {1-\th p^{-1} \over C} \,.
\ee
This shows that the cluster $\C(v_0)$ is infinite with positive probability.


The next step is to construct a flow $f$ of finite mean energy on $\C(v_0)$. 
For each $N \geq 1$ and every directed edge $e$ in $E_G$,  we define
\be \lab{eq:f}
f_N(e)= \int_{  \pas} p^{-N} 
\one_{\{\varphi_N \mbox{ \rm is open}\}}  \one_{\{ e \in \varphi_N\}}
 \, d\mu(\varphi) \,.
\ee
Then  $f_N$ is a flow on $\C(v_0) \cap B(v_0,r_N+1)$
   from $v_0$ to the complement of $B(v_0,r_N)$, i.e., for any vertex 
$v \in B(v_0,r_N)$ except $v_0$, the incoming flow to $v$ equals
the outgoing flow from $v$.
The {\em strength} of $f_N$  (the total outflow from $v_0$) is exactly $Z_N$.

Next, we estimate the expected energy of $f_N$:
\bea \lab{eq:energy}
\E \! \sum_{e \in E_G} f_N(e)^2 \! &=& \! \E \! \int_{ \pas} \!\! \int_{ \pas} p^{-2N} 
\one_{\{\varphi_N , 
\psi_N \mbox{ \rm \small are open}\}} \sum_{e \in E_G}
 \one_{\{ e \in \varphi_N\}} 
\one_{\{ e \in \psi_N\}}  \, d\mu(\varphi) \, d\mu(\psi) \\[1.4ex] \nonumber
 &=&\int_{ \pas} \!\! \int_{ \pas} |\varphi \cap \psi| \, p^{-|\varphi \cap \psi|}
 \, d\mu(\varphi) \, d\mu(\psi) \,.
\eea

Another application of exponential intersection tails
 allows us to bound the last integral by
\be \lab{eq:enbnd}
\sum_{n=0}^\infty C \th^n n p^{-n} \, ,
\ee
which is finite for $p>\th$.
For each edge $e$ of $G$, the sequence $\{f_N(e)\}$ 
is bounded in $L^2$, so it has a  weakly convergent subsequence.
Using the diagonal method, we can find a single increasing  sequence $\{N(k)\}_{k \ge 1}$
such that for every edge $e$, the sequence
 $f_{N(k)}(e)$ converges weakly as $k \to \infty$ to a limit, denoted $f(e)$.
Recalling that $Z_N$ is the strength of $f_N$,
we deduce that  $f_{N(k)}(e)$ converges weakly as $k \to \infty$ to a limit, denoted 
$Z_{\infty}$.
Since $r_N \to \infty$,  the limit function $f$ is a.s. a 
flow of strength $Z_\infty$ on $\C(v_0)$.
Exhausting $G$ by finite sets of edges, we conclude that 
the expected energy of $f$ is also bounded by (\ref{eq:enbnd}).
Thus 
$$
\P[\C(v_0) \mbox{ \rm is transient}] \geq \P[Z_\infty>0] >0 \, ,
$$
so the tail event $[\exists v : \C(v) \mbox{ \rm is transient}]$
must have probability $1$ by Kolmogorov's zero-one law.

Finally, we remove the assumption (\ref{eq:assump}). Any path tending to infinity
contains a self-avoiding path, obtained by ``loop-erasing''
(erasing cycles as they are created); see chapter 7 in Lawler (1991).
Thus we may indeed assume that $\mu$ is supported on self-avoiding paths.
 Since all paths
in $\pa$  tend to infinity, by Egorov's theorem there is a  closed subset $\pa'$
of $\pa$ on which this convergence is uniform, such that $\mu(\pa') >\mu(\pa)/2$.
Restricting $\mu$ to $\pa'$ and normalizing, we obtain a probability measure $\mu'$
on $\pa'$ that satisfies (\ref{eq:assump}) and (\ref{eq:eit}) with $4C$ in place of $C$,
so the proof given above applies.
$\Cox$

\noindent{\bf Remark.} Let $\pa_1 = \pa_1(G,v_0) \subset \pa$ denote the set
  of {\bf paths with unit speed}, i.e., those paths such that the $n^{th}$ vertex 
  is at distance $n$ from $v_0$, for every $n$.
In most applications of Proposition \ref{prop:trans}, the measure $\mu$ is supported
on $\pa_1$. When that is the case, the flows $f_N$ considered in the preceding
proof converge a.s.\  to a flow $f$, so there is no need to pass to subsequences.
Indeed, let $\B_N$ be the $\sif$ generated by the status (open or closed)
of all edges on paths  $\varphi_N$ with $\varphi \in \pa$.
 It is easy to check that $\, \{Z_N\}_{N \geq 1}$ is
an $L^2$ martingale adapted to the filtration $\{\B_N\}_{N \geq 1}$.
Therefore $\{Z_N\}$ converges a.s.\ and in $L^2$ to
a mean 1 random variable $Z_\infty$. Moreover, for each edge $e$ of $G$, 
the sequence $\{f_N(e)\}$ is a $\{\B_N\}$-martingale which 
converges a.s.\ and in $L^2$ to a nonnegative random variable $f(e)$. 

\begin{cor}[Grimmett, Kesten and Zhang (1993)] \lab{cor:GKZ}
Consider ordinary bond percolation with parameter $p$ on $\Z^d$, where $d \geq 3$.
For all $p>p_c$, the unique infinite cluster is a.s.\ transient.
\end{cor}

\noindent{\sc Proof:}\,
Transience of the unique infinite cluster is a tail event, so it has probability 0 or 1.
Since $\Z^d$ admits random paths with EIT for $d \geq 3$ 
(see Theorem \ref{thm:3d} and the discussion preceding it),
it follows from Proposition \ref{prop:trans} that the infinite cluster is transient 
if $p$ is close enough to 1. As remarked before, this conclusion also
applies to site percolation. 

Recall that a set of graphs ${\bf B}$ is called {\bf increasing} if for any 
graph $G$ that contains a subgraph in  ${\bf B}$, necessarily
$G$ must also be in  ${\bf B}$.

Consider now bond percolation with {\em any} parameter $p>p_c$
in $\Z^d$. Following Pisztora (1996), call an open cluster $\C$ contained in some cube $Q$
a {\em crossing cluster\/} for $Q$ if for all $d$ directions there is an open path
contained in $\C$ joining the left face of $Q$ to the right face. 
For each $v$ in the lattice $N \Z^d$, denote by $Q_N(v)$ the cube of side-length
$5N/4$ centered at $v$, and let
$A_p(N)$ be the set of $v \in N \Z^d$ with the following property: \, \em
The cube $Q_N(v)$
contains a crossing cluster $\C$ such that any open cluster in $Q_N(v)$
of diameter greater than $N/10$ is connected to $\C$ by an open path
in $Q_N(v)$.
 \rm

Proposition 2.1 in Antal and Pisztora (1996), which relies on the
work of Grimmett and Marstrand (1990), implies that 
$A_p(N)$ stochastically dominates site percolation  with parameter
$p^*(N)$ on the stretched lattice $N \Z^d$, where $p^*(N) \to  1$ as 
$N \to \infty$.
(Related renormalization arguments can be found in 
 Kesten and Zhang (1990) and Pisztora (1996);
 General results on domination by product measures where obtained by
 Liggett, Schonmann and Stacey (1996)).
This domination means that for any increasing Borel set 
of graphs ${\bf B}$, the probability that the subgraph 
of open sites under
independent site percolation  with parameter
$p^*(N)$ lies in ${\bf B}$,
is at most $\P[A_p(N) \in {\bf B}]$. 
If $N$ is sufficiently large, then the infinite cluster determined by the
site percolation with parameter $p^*(N)$ on the lattice
$N \Z^d$, is transient a.s.
By Rayleigh's monotonicity principle (see Doyle and Snell 1984),
the set of subgraphs of $N \Z^d$ that have a 
transient connected component is increasing, so $A_p(N)$ has a transient 
component $\widehat{A}_p(N)$ with probability 1.


Recall from  Doyle and Snell (1984) that the ``$k$-fuzz'' of
a graph $\Gamma=(V,E)$ is the graph  $\Gamma_k=(V,E_k)$
where the vertices $v,w \in V$ are connected by an edge in $E_k$
iff there is a path of length at most $k$ between them in $\Gamma$.
The $k$-fuzz $\Gamma_k$ is transient iff $\Gamma$ is transient
(See Section 8.4 in  Doyle and Snell (1984), or Lemma 7.5 in Soardi (1994).)
By assigning to each  $v \in N \Z^d$ a different vertex $F(v)$ in  the
intersection  $\C_p \cap Q_N(v)$, we see that
$\widehat{A}_p(N)$ is isomorphic to a subgraph of the
$2(5N/4)^d$-fuzz of the infinite cluster $\C_p$ in the original lattice.
It follows that $\C_p$ is also transient a.s.
$\Cox$

\section{A summable predictability profile yields EIT} 
The following lemma will imply Theorem \ref{thm:3d}.
\begin{lem} \lab{lem:pred}
Let $\{\Ga_n\}$ be a sequence of random variables taking values in a
countable set $V$.
If the predictability profile (defined in (\ref{eq:defpre})) of $\Ga$ satisfies
$ \sum_{k=1}^\infty \pr_\Gamma (k) < \infty $,
then there exist $C< \infty$ and $0<\th<1$ such that
for any sequence $\{v_n\}_{n \geq 0}$ in $V$ and all $\ell \geq 1$,
\be \lab{eq:sump}
\P [\#\{n \geq 0\, : \,  \Ga_n=v_n\} \geq \ell] \leq  C \th^\ell  \,.
\ee

\end{lem}

\noindent{\sc Proof:}
Choose $m$ large enough so that $ \sum_{k=1}^\infty \pr_\Gamma (km)= \beta<1$,
whence for any sequence $\{v_n\}_{n \geq 0}$ ,
$$
\P\Big[\exists k \geq 1 \, : \, \Ga_{n+km}= v_{n+km} \, \Big| \,  
\Gamma_0,\ldots,\Gamma_n\Big]  \leq \beta 
\quad \mbox{ \rm for all }  \: n \geq 0 \,.
$$
It follows by induction on $r \geq 1$ that for all 
$j \in \{0,1, \ldots, m-1 \}$,
\be \lab{eq:modm}
\P [\#\{k \geq 1 \, : \,  \Ga_{j+km}=v_{j+km}\} \geq r] \leq  \beta^{r} \, . 
\ee
If $\#\{n \geq 0 \, : \,  \Ga_n=v_n\} \geq \ell$ then there must be some
$j \in \{0,1, \ldots, m-1 \}$ such that 
$$
\#\{k \geq 1 \, : \,  \Ga_{j+km}=v_{j+km}\} \geq \ell/m-1 \, .
$$
Thus the inequality (\ref{eq:sump}), with $\th= \beta^{1/m}$ and 
$C=m\beta^{-1}$, follows from (\ref{eq:modm}). 
$\Cox$

\noindent{\bf Proof of Theorem \ref{thm:3d}.}\,
Let $\{S_n\}_{n \geq 0}$  be a nearest-neighbor process on $\Z$
starting from $S_0=0$, that satisfies 
(\ref{eq:smooth}) for some $\al>1/2$ and $C_\al < \infty$. 
Denote $W_n=(n+S_n)/2$ and suppose that the processes $\{W_n\}$
and  $\{W^{\sharp}_n\}$
are independent of each other and have the same distribution.
Clearly $\{W_n\}$ and $\{S_n\}$ have the same predictability profile.
We claim that   \em the random oriented path  
$$
\{\Gamma_n\}_{n \geq 0} = \left\{ \, \left(W_{\lfloor n/2 \rfloor},  
   \, W^{\sharp}_{\lceil n/2 \rceil}, \, 
   n-W_{\lfloor n/2 \rfloor}-W^{\sharp}_{\lceil n/2 \rceil}\right) \,\right\}_{n \geq 0}
$$
in $\, \Z^3\,$ has exponential intersection tails.

\rm First, observe that $\{\Gamma_n\}$ is indeed an oriented path,
i.e., $\Gamma_{n+1}-\Gamma_n$ is one of the three vectors
$(1,0,0), (0,1,0), (0,0,1)$ for every $n$.
Second,  $\pr_{\Ga}(k)= \pr_{S}(\lfloor k/2 \rfloor)^2
\leq C_\al^2 \lfloor k/2 \rfloor^{-2\al}$ is summable in $k$,
so $\{\Gamma_n\}$ satisfies (\ref{eq:sump}) for some $C< \infty$ and $0<\th<1$.
Denote the distribution of  $\{\Gamma_n\}$ by $\mu$, and let
$\{\Gamma^*_n\}$ be an independent copy of $\{\Gamma_n\}$.
Integrating (\ref{eq:sump}) with respect to $\mu$, we get
$$
\forall \ell \quad 
\mu \times \mu \Big[\#\{n \geq 0\, : \,  \Ga_n=\Ga^*_n\} \geq \ell\Big] \leq 
 C \th^\ell  \,.
$$
For two oriented paths $\Ga$ and $\Ga'$ in $\Z^d$, the number of edges in 
common is at most the collision number $\#\{n \geq 0\, : \,  \Ga_n=\Ga^*_n\}$,
 and hence $\mu$ has EIT($\th)$.
$\Cox$

\section{Summing boundary spins yields an unpredictable path}
In this section we prove Theorem \ref{thm:smooth}. 
The engine of the proof is Lemma \ref{lem:recur} concerning 
the distribution of 
the population in a two-type branching process.
 Let $ \ell \geq 2$ and  $r \geq 1$
 be integers and write $b=\ell+r$.
Denote by 
$\T_b$  the  infinite rooted tree where each vertex has exactly $b$ children.
Consider the following 
 labeling $\{\sigma(v)\}$ of the vertices of $\T_b$ by $\pm 1$ valued random variables,
called {\em spins} because of the analogy with the Ising model
(cf. Moore and Snell 1979).
Let $\sigma({\rm root})=1$. 
For any vertex $v$ of $\T_b$ with children $w_1,\ldots, w_b$, 
assign the first $\ell$ children the same spin as their
parent:  $\sigma(w_j)=\sigma(v)$ for $j=1,\ldots,\ell$,
 and assign the other  $r$ children i.i.d.  spins 
$ \sigma(w_{\ell+1}) , \ldots, \sigma(w_b)$
that take the values $\pm 1$ with equal probability, and are independent of 
$\sigma(v)$.
 As $N$ varies, the population vectors $(Z_N^+,Z_N^-)$
which count the number of spins of each type at level $N$ of $\T_b$,
form a two-type branching process with mean offspring matrix 
$$
\left( \ba{lr}  \ell+r/2 & r/2 \\
                  r/2  & \ell+r/2 \ea \right) \, \, .
$$
(See Athreya and Ney (1972) for background on branching processes).
The Perron eigenvalue of this matrix is $b$, but we are interested in the scalar product
of the population vectors with the eigenvector $(1,-1)$ which corresponds
to the second eigenvalue $\ell$ of the mean offspring matrix.
\begin{lem} \lab{lem:recur}
  The sum of all spins at level $N$
\be \lab{eq:summ}
Y_N= \sum_{|v|=N} \sigma(v) = Z_N^+-Z_N^-
\ee
satisfies the inequality
\be \lab{eq:bnd}
\P[Y_N=x] \leq C \ell^{-N}
\ee
for all $N \geq 1$ and all integer $x$, where the constant $C$ depends only on $\ell$
and $r$.
\end{lem} 
{\bf Remark:} A closer examination of the proof below  shows that
the inequality (\ref{eq:bnd}) is sharp
iff $\ell> \sqrt{b}$. The significance of the condition $\ell> \sqrt{b}$ 
is explained by Kesten and Stigum (1966) in a more general setting. If 
$\ell < \sqrt{b}$
then the distributions of $Y_N b^{-N/2}$ converge to a normal law.

\noindent{\bf Proof of Lemma \ref{lem:recur}:}\, 
By decomposing the sum in the definition of $Y_{N+1}$
into $b$ parts according to the first level of $\T_b$, we see that
\be \lab{eq:rec}
Y_{N+1}= \sum_{j=1}^\ell  Y_N^{(j)} + \sum_{j=\ell+1}^b \sigma(w_j) Y_N^{(j)}  \, ,
\ee
where $\{\sigma(w_j)\}_{j=\ell+1}^b$ are $\, r \:$ i.i.d.\ uniform spins,
and $\{Y_N^{(j)}\}_{j=1}^b$ are i.i.d.\ variables
with the distribution of $Y_N$, that are independent of these spins.
Consequently, the characteristic functions
\be \lab{eq:char}
\wh{Y}_N (\lambda)= \E(e^{i \lambda Y_N})
\ee
satisfy the recursion
\be \lab{eq:char2}
\wh{Y}_{N+1} (\lambda) = 
\wh{Y}_N (\lambda)^\ell  \Big( \frac{\wh{Y}_N (\lambda)+\wh{Y}_N (-\lambda)}{2} \Big)^r =
 \wh{Y}_N (\lambda)^\ell  ( \Re \wh{Y}_N (\lambda))^r \,,
\ee
where $\Re$ denotes real part. Using the polar representation 
$\wh{Y}_N = |\wh{Y}_N| e^{i\gamma_N(\la)}$,
the last equation implies that   
$\gamma_{N+1}(\la) \equiv \ell \gamma_N(\la) \, \bmod \pi$.
(Note that $\Re \wh{Y}_N (\lambda)$ may be negative.)
By definition $Y_0=1$, and therefore $\wh{Y}_0(\lambda) = e^{i\la}$, 
so $\gamma_0(\la)=\la$. Consequently,
$\gamma_N(\la) \equiv\ell^N \la  \bmod \pi$ for all $N$. 
Taking absolute values in  (\ref{eq:char2}) yields
\be \lab{eq:char3}
|\wh{Y}_{N+1} (\lambda)| = |\wh{Y}_N (\lambda)|^b \cdot|\cos (\ell^N \la) |^r \,.
\ee
By induction on $N$, we obtain
\be \lab{eq:char4}
\forall N \geq 0 \quad \:
|\wh{Y}_{N} (\lambda)| = \prod_{k=1}^{N} |\cos (\ell^{N-k} \la) |^{rb^{k-1}} \,.
\ee
Since $|\wh{Y}_N(\cdot)|$ is an even function with period $\pi$, 
changing variables $t=\ell^N \la$ gives
\be \lab{eq:char5}
\int_{-\pi}^\pi |\wh{Y}_{N} (\lambda)| \, d\la = 
4 \int_{0}^{\pi/2} |\wh{Y}_{N} (\lambda)| \, d\la =
4  \ell^{-N} \int_{0}^{\frac{\ell^N \pi}{2} } 
 \prod_{k=1}^{N} |\cos (\ell^{-k} t) |^{rb^{k-1}} \, dt
\ee
Denote $\xi= \cos \frac{\pi}{2 \ell}$. 
For $t \in [\ell^{k-1}\pi/2 ,\ell^k \pi/2]$, the $k$th factor in the rightmost integrand
in (\ref{eq:char5}) is bounded by $\xi^{rb^{k-1}}$, 
and the other factors are
at most $1$. Consequently,
\be \lab{eq:char6}
\int_{-\pi}^\pi |\wh{Y}_{N} (\lambda)| \, d\la \leq 
4\ell^{-N}(\frac{\pi}{2}+  \sum_{k=1}^{N} \frac{\ell^k \pi}{2} 
\xi^{rb^{k-1}}) \leq 
C \ell^{-N}
\ee
where $C$ depends only on $\ell$ and $r$, since the sum 
$\sum_{k=1}^{\infty} \ell^k  \xi^{rb^{k-1}}$ converges.
Finally, Fourier inversion and (\ref{eq:char6}) imply that for any integer $x$,
$$
\P[Y_N=x] = \frac{1}{2\pi}
 \int_{-\pi}^\pi \wh{Y}_{N} (\lambda) e^{-i\la x} \, d\la \leq
\int_{-\pi}^\pi |\wh{Y}_{N} (\lambda)| \, d\la \leq  C \ell^{-N} \, .
$$
$\Cox$

\noindent{\bf Proof of Theorem \ref{thm:smooth}:}\,\newline  {\bf (a)}
Given integers $\ell \geq 2$ and $r \geq 1$, let $b=\ell+r$ as above.
Embed the regular tree $\T_b$ in the upper half-plane, 
with the root on the real line and the children of every vertex arranged
    from left to right above it. Label the
vertices of $\T_b$ by $\pm 1$ valued spins $\{\sigma(v)\}$ 
as described at the beginning of this section.
 Let $M>1$ and suppose that
$b^N \geq M$. Let $\{v_j\}_{j=1}^{b^N}$ be  the vertices at level $N$
of $\T_b$, enumerated from left to right. For $m \leq M$,
denote $S_m= \sum_{j=1}^m \sigma(v_j)$ and observe that the joint 
distribution of the $M$ random variables $\{S_m\}_{m=1}^M$ does not depend on
 $N$. Using Kolmogorov's Consistency Theorem, we obtain an infinite process
 $\{S_m\}_{m=1}^\infty$. We claim that the predictability
profile of this process satisfies
\be \lab{eq:smooth2}
\pr_S(k) \leq (2b)^\al C k^{-\al} \quad \mbox{ \rm for all } k \geq 1 \, ,
\ee
where $\al=\frac{\log \ell}{\log b}$, and $C=C(\ell,r) \geq 1$ is given in (\ref{eq:bnd}).
Since we can take $\ell$ arbitrarily large and $r=1$,
establishing (\ref{eq:smooth2}) for $k>2$ will suffice to prove the theorem.

Given $n \geq 0$ and $k \geq 1$, choose $N$ such that $b^N \geq n+k$,
so the random variables $\{S_j\}_{j=0}^{n+k}$ may be obtained by 
summing spins along level $N$ of $\T_b$.
There is a unique $h \geq 0$ such that
\be \lab{eq:defh}
   2b^h \leq k  <2b^{h+1} \,.
\ee
For any vertex $v$, denote by $|v|$ its level in $\T_b$, and for $i \geq 1$ 
let  $D_i(v)$ be the set  of 
its $b^i\,$ descendants at level $|v|+i$.
By (\ref{eq:defh}), there exists at least one vertex $v$ at level $N-h$ of 
$\T_b$, such that
$D_h(v)$ is contained in  $\{v_{n+1}, v_{n+2}\ldots,v_{n+k}\}$
in the left-to-right enumeration of level $N$. 
Denote by $D(v) = \bigcup_{i=1}^\infty D_i(v)$
the set of all descendants of $v$, 
and by $\F_v^*$ the sigma-field generated by all the
spins $\{\sigma(w) : w \notin D(v)\}$. The random variable
$$
\tilde{Y}_h(v) =\sigma(v) \sum_{w \in D_h(v)} \sigma(w)
$$
is independent  of $\F_v^*$, and  has the same distribution as the variable 
$Y_h$ defined by (\ref{eq:summ}).
Clearly, we can write 
$$
S_{n+k} = \sigma(v) \tilde{Y}_h(v) + S^*_{n+k} \, ,
$$
where $S^*_{n+k}$, the sum of $n+k-b^h$ spins labeling vertices not in $D_h(v)$, is  $\F_v^*$-measurable.
Consequently, for any integer $x$,
\be \lab{eq:almost}
\P[S_{n+k}=x \mid \F_v^*] = \P[\tilde{Y}_h(v)= \sigma(v)(x-S^*_{n+k}) \mid \F_v^*] \leq C \ell^{-h} \,,
\ee
by Lemma \ref{lem:recur}.
The definitions  of $\alpha$ and $h$ imply that
$\ell^{-h} = b^{-h \al}$ and $b^{h} > k/2b$, so we infer from
(\ref{eq:almost}) that
$$
\forall x \in \Z \quad \P[S_{n+k}=x \mid \F_v^*] < (2b)^\al C k^{-\al} \,.
$$
Since $S_0, S_1, \ldots S_n$ are $ \F_v^*$-measurable, this yields 
(\ref{eq:smooth2})
and completes the proof of part (a) of the theorem. \newline
{\bf (b)}  The property (\ref{eq:smooth}) is stable under shifts, mixtures, weak limits, and passing to ergodic components,
 so it is possible to obtain the desired stationary 
process as an ergodic component of a weak limit point of the averages
$\frac{1}{n}(S+\Theta S+ \cdots +\Theta^{n-1} S)$, where $\Theta$ is the left
shift. 

We now describe such a process more explicitly, by modifying
the construction in part (a).
Let $\sigma({\rm root})$ be a uniform random spin; define the other spins
  from it as in part (a).
Given $N>1$,
choose $U$ uniformly in $\{1, \ldots, b^N\}$ and define
$$
\tilde{S}_n= \sum_{j=U}^{U+n-1} \sigma(v_j)
\quad \mbox{ \rm for } n \leq b^N-U+1 \,,
$$
where  $\{v_j\}_{j=1}^{b^N}$ is the left-to-right
enumeration of level $N$ of $\T_b$.
To extend the sequence $\tilde{S}$ further,
we consider the root of $\T_b$ as the $J$th child $w_J$ of a new vertex $\rho$,
where $J$ is chosen uniformly in $\{1, \ldots,b\}$. If $J \leq \ell$,
let $\sigma(\rho)= \sigma(w_J)$,
and if $J >\ell$ let $\sigma(\rho)$ be a uniform random spin, 
independent of the spins on the original
tree. We can view the original tree $\T_b$ as a subtree of
a  new $b$-tree $\tilde{\T}_b$ rooted at $\rho$.
Since $(J-1)b^N+U$ is uniformly distributed in $\{1, \ldots, b^{N+1}\}$,
the vertex $v_{_U}$ is uniformly distributed in level $N+1$
of $\tilde{\T}_b$.
Repeating this re-rooting procedure and enlarging $N$  as needed,
yields the desired process $\{\tilde{S}_j\}_{j=1}^\infty$.
The proof given in part (a) also shows that this process has the
unpredictability property (\ref{eq:smooth}).
Stationarity and ergodicity of the increments can be derived from the 
invariance and ergodicity of the  Haar measure on the 
$b$-adic integers under the operation of adding 1; we omit the details.
$\Cox$
%

\section{Concluding remarks and questions}
\begin{enumerate} 
\item Consider the following three properties that an infinite connected
 graph $G$ may have:
\begin{description}
\item{{\bf (i)}} $G$ admits random paths with EIT.
\item{{\bf (ii)}}  There exists $p < 1$ such that
 simple random walk is transient  on a percolation cluster of $G$ for
bond percolation with  parameter  $p$.
\item{{\bf (iii)}}  A random walk in 
random environment on $G$ defined by i.i.d.\ resistances
with any common distribution is almost surely transient.
\end{description}
In Pemantle and Peres (1996) it is shown that
properties (ii) and (iii) are equivalent.
Proposition \ref{prop:trans} of the present paper shows that (i) implies (ii);
does (ii) imply (i)? \, \newline
Note that there exist transient trees of polynomial growth
(see, e.g., Lyons 1990),
and these cannot admit random paths with EIT since they have $p_c=1$.
\item Does $\Z^d$ with $d \geq 3$ admit random paths with EIT$(\th)$ 
for {\em all}
$\th>p_c$? \newline
(This question was suggested to us by Rick Durrett.)
A similar question can be asked for other graphs in place of
$\Z^d$ , e.g., for transient Cayley graphs.
A positive answer to this question when the graph in question is a {\em tree}
follows from the work of Lyons (1990). Indeed, a flow from the root of the tree
can be identified with a measure on paths, and the energy of the flow $\mu$ 
in the kernel $p^{-|x \wedge y|}$ can be
identified with an exponential moment of the number
of intersections of two paths chosen independently according to $\mu$.
\item Lyons (1995) finds a tree with $p_c<1$ contained 
as a subgraph in the Cayley graph of
any group of exponential growth. It follows that such Cayley graphs
admit paths with EIT.
\item \label{it:wedges} It is easy to adapt the proofs of Theorem \ref{thm:3d} and 
Corollary \ref{cor:GKZ} to show that for any $\epsilon>0$,
the cone $\{ (x,y,z) \in \Z^3 \, :\, |z| \leq  \epsilon |x| \}$ in $\Z^3$
admits random paths with EIT, and contains a transient percolation cluster
for all $p>p_c$. Does the subgraph
 $\{ (x,y,z) \in \Z^3 \, :\, |z| \leq   |x|^\epsilon \}$ share these 
properties?\, (This subgraph is sometimes viewed as a model for
a ``$2+\epsilon$~dimensional lattice''.)
\item Does oriented percolation in $\Z^d$ admit transient infinite clusters
  $\C(v)$ for all parameters $p>p_c^{\mbox{\rm \small  or }}$? \,
The challenge here is to adapt the renormalization argument 
used in the proof of Corollary \ref{cor:GKZ} to the oriented setting. 
\item Consider the stationary processes $\{\tilde{S}_n\}$
constructed at the end of the previous section, and let $\al=\ell/b$
with $\sqrt{b}< \ell <b$.
 Do the rescaled step functions $t \mapsto n^{-\al} \tilde{S}_{\lfloor nt \rfloor}$
on $[0,1]$ converge in law? \,  Is the limit a Gaussian process?  \, 
It is easily verified that $\E|\tilde{S}_n-\tilde{S}_m|^2 \asymp |n-m|^{2\al}$,
which is reminiscent of fractional Brownian motion.
 The proof of Lemma \ref{lem:recur} implies that
 $b^{-\al N} Y_N$ converges in law to a (non-Gaussian)
distribution with characteristic function
$ s \mapsto e^{is} \prod_{k=1}^{\infty} [\cos (\ell^{-k} s) ]^{rb^{k-1}}$ \,
(Recall that $\ell>\sqrt{b}$).
\item \label{it:erg} How fast can the predictability profile 
(\ref{eq:defpre}) of a nearest-neighbor
process on $\Z$ decay? \, 
By Theorem \ref{thm:smooth}, a decay rate of $O(k^{-\alpha})$  
is possible for any $\al <1$.
 On the other hand, a decay rate of $O(1/k)$ is impossible.
Indeed, if there exists a  nearest-neighbor process
with predictability profile bounded by $\{g(k)\}$, then there exists
such a process with stationary ergodic increments;
then $g(k)= O(1/k)$ is ruled out by the ergodic theorem.

\item Among nearest-neighbor processes $\{S_n\}$ 
on $\Z$, clearly simple random walk
has the most unpredictable {\em increments}, in any conceivable sense.
Heuristically, there is a tradeoff here:
when the increments are very unpredictable (e.g., their predictability
profile tends rapidly to $1/2$), cancellations dominate,
and the partial sums becomes more predictable.
Our construction in Section 4 sacrificed the independence 
of the increments, to make their partial sums less predictable.
It would be quite interesting to establish a precise
quantitative form of this tradeoff. 

\item Is there a  construction of a measure on paths in $\Z^3$
with exponential intersection tails, which is simpler 
than that given in Sections 3 and 4?
\end{enumerate}

\subsection{Recent developments}

 After a previous version of the present paper was circulated, 
some of the problems raised above were solved,
and several further extensions of the fundamental transience theorem of 
Grimmett, Kesten and
Zhang (1993) were obtained.
\begin{itemize}
\item 
H\"{a}ggstr\"{o}m and Mossel (1998) 
constructed processes
with predictability profiles bounded by $C/[k(f(k)]$,
for any decreasing $f$ such that $\sum_j f(2^j) < \infty$.
They gave two different constructions, one based on the Ising model
on trees, and the other via a random walk with a random drift that varies
in time. H\"{a}ggstr\"{o}m and Mossel also answered affirmatively Question \ref{it:wedges}
 above, by
constructing paths with exponential intersection tails in ``$2+\epsilon$'' dimensions.
Remarkably, they were able to show that for a class of
 ``trumpet-shaped'' subgraphs $G$ of $\Z^3$, transience of $G$ implies 
a.s.\ transience of an infinite  percolation cluster  in $G$ for any $p>p_c$.

\item In a brief but striking paper, Hoffman (1998) improved the bounds in Remark 
\ref{it:erg}  above, and showed that the constructions of H\"{a}ggstr\"{o}m and Mossel 
are optimal. Specifically, he used a novel renormalization argument to
prove  that if $f$ satisfies $\sum_j f(2^j) = \infty$, then 
there is no nearest--neighbor 
process on $\Z$ with predictability profile bounded by $f$.  

\item Hiemer (1998) proved a renormalization theorem for oriented
percolation, that allowed him to extend our result on transience of oriented
percolation clusters in $\Z^d$ for  $d \ge 3$, from the case of high $p$
to the whole supercritical phase $p>p_c^{\mbox{\rm \small  or }}$.

\item  
Consider supercritical percolation on $\Z^d$ for $d \geq 3$.
The transience result of Grimmett, Kesten and Zhang (1993) 
 is equivalent to the existence of a nonzero flow $f$ on the infinite 
cluster such that the $2$--energy $\sum_{e}f(e)^2$ is finite.
Using the method of unpredictable paths, 
Levin and Peres (1998) sharpened this result, and showed that 
the infinite cluster supports a nonzero flow $f$ such that the  
$q$--energy $\sum_{e}|f(e)|^{q}$  is finite for all $q>d/(d-1)$.   
Thus the infinite cluster has the same ``parabolic index''
as the whole lattice. 
(See the last chapter of Soardi (1994) for the definition.)

\end{itemize}
\noindent{\bf Acknowledgement:} We are grateful to Geoffrey Grimmett
for a helpful remark  on renormalization. We also thank
Chris Hoffman, David Levin, Elhanan Mossel and the referee
for their useful comments on the presentation.

\renewcommand{\baselinestretch}{1.0}\large


\normalsize
\begin{thebibliography}{YMNX}
\bibitem{AP}
P. Antal and A. Pisztora (1996). On the chemical distance in supercritical Bernoulli 
 percolation. {\em Ann. Probab.\  } {\bf 24} {1036--1048}.

\bibitem{AN}
K. Athreya and P. Ney (1972).  {\em Branching Processes.} 
 Springer-Verlag, New York.

\bibitem{CD}
T. Cox and R. Durrett (1983). Oriented percolation in dimensions $d \geq 4$:
bounds and asymptotic formulas. {\em Math.\ Proc.\ Camb.\ Phil.\ Soc.\ }
{\bf 93}, 151--162.

\bibitem{DS}
P. G. Doyle and E. J. Snell (1984).
{\em Random walks and electrical networks.}
 Carus Math.\ Monographs {\bf 22}, Math.\ Assoc.\ Amer., Washington, D. C.

\bibitem{GKZ}
G. R. Grimmett, H. Kesten and Y. Zhang (1993).
Random walk on the infinite cluster of the percolation model.
{\em Probab.\ Th.\ Rel.\ Fields} {\bf 96}, 33--44.

\bibitem{GM} G. R. Grimmett and J. M. Marstrand (1990).
The supercritical phase of percolation is well behaved.
{\em Proc.\ Royal Soc. London Ser. A}  {\bf 430}, 439--457.

\bibitem{HM} O. H\"{a}ggstr\"{o}m and E. Mossel (1998).
  Nearest-neighbor walks with low predictability profile and
  percolation in $2+\epsilon$ dimensions. {\em Ann.\ Probab.}, to appear. 

\bibitem{Hi} P. Hiemer (1998).  Dynamical renormalisation in oriented
  percolation.  {\em Preprint.}

\bibitem{Ho} C. Hoffman (1998).  Unpredictable nearest neighbor
processes. {\em Preprint.}


\bibitem{KS} H. Kesten and B. P. Stigum (1966).
Additional limit theorems for indecomposable multidimensional
Galton-Watson processes. {\em Ann.\ Math.\ Stat.\ } {\bf 37},
1463--1481.


\bibitem{KZ} H. Kesten and Y. Zhang (1990). The probability 
 of a large finite cluster in supercritical Bernoulli percolation.
{\em Ann.\ Probab.\  } {\bf 18}, 537--555.

\bibitem{La}
G. Lawler (1991).
{\em Intersections of random walks.} Birkh\"{a}user, Boston.

\bibitem{LP} D. Levin and Y. Peres (1998). 
Energy and Cutsets in Infinite Percolation Clusters. 
 {\em  Proceedings of the Cortona Workshop on Random Walks
and Discrete Potential Theory}, M. Picardello and W. Woess (editors),
to appear.

\bibitem{LSS} T. M. Liggett, R. H. Schonmann and A. M. Stacey (1996).
Domination by product measures. {\em Ann.\ Probab.} {\bf 24},  1711--1726.

\bibitem{Ly90}
R. Lyons (1990).  Random walks and percolation on trees.
{\em Ann.\ Probab.\  } {\bf 18} 931--958.

\bibitem{Ly95}
R. Lyons (1995). Random walks and the growth of groups.
{\em C. R. Acad.\ Sci.\ Paris} {\bf 320}, 1361--1366. 


\bibitem{MS}
T. Moore and J. L. Snell (1979).  A branching process showing a phase transition.
{\em J. Appl.\ Probab.} {\bf 16}, 252--260. 

\bibitem{PP} R. Pemantle and Y. Peres (1996). On which graphs are all 
random walks
in random environments transient? {\em Random Discrete
Structures}, IMA Volume 76 (1996), D. Aldous and R. Pemantle (Editors), 
  Springer-Verlag.

\bibitem{Pi}
A. Pisztora (1996). Surface order large deviations for Ising, Potts and percolation
models. {\em Probab.\ Th.\ Rel.\ Fields} {\bf 104}, 427--466.

\bibitem{So}
P. M. Soardi (1994). {\em Potential Theory on Infinite Networks}. 
Springer LNM, Berlin.

\end{thebibliography}
\end{document}